\documentclass[12pt]{amsart}

\title{Long monotone paths on simple 4-polytopes}

\author{Julian Pfeifle}

\address{Institut de Matemàtica, Universitat de Barcelona, Gran Via de
  les Corts Catalanes 585, E-08007 Barcelona, Spain }

\date{February 16, 2004}
\email{julian@imub.ub.es}

\thanks{The author was financed by the DFG Graduiertenkolleg
  \emph{Combinatorics, Geometry, and Computation} (GRK 588-2), the GIF
  project \emph{Combinatorics of Polytopes in Euclidean Spaces}
  (I-624-35.6/1999), and post-doctoral fellowships from MSRI and
  Institut de Matemàtica de la Universitat de Barcelona.    }

\usepackage{amsmath}        
\usepackage{amssymb}
\usepackage{amsthm}
\usepackage{graphicx}
\usepackage{color}
\usepackage{psfrag}         
\usepackage{bbm}            
\usepackage{lscape}          
\usepackage[small,sc]{caption2}
\usepackage{alltt}
\usepackage[noend]{algorithmic}
\usepackage{algorithm}      
\usepackage{tabularx}       
\usepackage{ae,aecompl}     
\usepackage[newitem,newenum]{paralist}
\usepackage{a4wide}

\renewcommand{\arraystretch}{1.2}
\newcommand{\plabel}[1]{({M}#1)}

\graphicspath{{pictures/}}   
\makeatletter
\def\input@path{{pictures/}}    
\makeatother

\newtheorem{theorem}{Theorem}[section]

\newtheorem*{theorem*}{Theorem}
\newtheorem*{maintheorem}{Main Theorem}

\theoremstyle{definition}
\newtheorem{conjecture}[theorem]{Conjecture}

\newtheorem{convention}[theorem]{Convention}

\newtheorem{definition}[theorem]{Definition}

\newtheorem{lemma}[theorem]{Lemma}
\newtheorem{notation}[theorem]{Notation}
\newtheorem{observation}[theorem]{Observation}
\newtheorem{proposition}[theorem]{Proposition}
\newtheorem{remark}[theorem]{Remark}

\newtheorem*{assumption*}{Assumption}

\newtheorem*{corollary*}{Corollary}
\newtheorem*{example*}{Example}
\newtheorem*{lemma*}{Lemma}

\newtheorem*{proposition*}{Proposition}
\newtheorem*{remark*}{Remark}

\renewcommand{\Pr}{\mathbbm{P}}
\newcommand{\R}{\mathbbm{R}}

\newcommand{\N}{\mathbbm{N}}

\DeclareMathOperator{\aff}{aff}
\DeclareMathOperator{\conv}{conv}

\DeclareMathOperator{\relint}{relint}

\DeclareMathOperator{\sk}{sk}

\DeclareMathOperator{\vertices}{vert}

\newcommand{\pos}{^{\ge0}}

\renewcommand{\phi}{\varphi}

\newcommand{\e}{\varepsilon}

\newcommand{\F}{{\mathcal F}}

\newcommand{\cH}{{\mathcal H}}

\newcommand{\calO}{{\mathcal O}}

\renewcommand{\a}{{\boldsymbol{a}}}

\newcommand{\be}{{\boldsymbol{e}}}

\newcommand{\n}{{\boldsymbol{n}}}

\newcommand{\q}{{\boldsymbol{q}}}

\newcommand{\s}{{\boldsymbol{s}}}

\newcommand{\bu}{{\boldsymbol{u}}}
\newcommand{\bv}{{\boldsymbol{v}}}

\newcommand{\x}{{\boldsymbol{x}}}

\newcommand{\balpha}{{\boldsymbol{\alpha}}}
\newcommand{\bbeta}{{\boldsymbol{\beta}}}
\newcommand{\bomega}{{\boldsymbol{\omega}}}

\newcommand{\bt}{\tilde{b}}
\newcommand{\vt}{\tilde{v}}
\newcommand{\wt}{\tilde{w}}
\newcommand{\Qt}{{\widetilde{Q}}}
\newcommand{\E}{\widetilde{E}}

\newcommand{\pit}{\widetilde{\pi}}

\newcommand{\vb}{\overline{v}}
\newcommand{\wb}{\overline{w}}

\begin{document}

\def\currentvolume{} \def\currentissue{} \pagespan{1}{60} \PII{}
\copyrightinfo{}{} \keywords{} 

\begin{abstract}
  The \emph{Monotone Upper Bound Problem} (Klee, 1965) asks
  if the number~$M(d,n)$ of vertices in a monotone path along edges of
  a $d$-dimensional polytope with $n$~facets can be as large as conceivably
  possible: Is $M(d,n)=M_{\rm ubt}(d,n)$, the maximal number of
  vertices that a $d$-polytope with $n$~facets can have according to
  the Upper Bound Theorem?
  
  We show that in dimension $d=4$, the answer is ``yes'', despite the
  fact that it is ``no'' if we restrict ourselves to the
  dual-to-cyclic polytopes. For each $n\ge5$, we exhibit a realization
  of a polar-to-neighborly $4$-dimensional polytope with $n$~facets
  and a Hamilton path through its vertices that is monotone with
  respect to a linear objective function.
  
  This constrasts an earlier result, by which no polar-to-neighborly
  $6$-dimensional polytope with $9$~facets admits a monotone Hamilton
  path. 
\end{abstract}

\maketitle

\section{Introduction}

While investigating the complexity of the simplex algorithm for linear
programming, Klee \cite{Klee65} in 1965 posed the \emph{Monotone Upper
  Bound Problem}: For $n>d\ge2$, he asked for the maximal number
$M(d,n)$ of vertices of a $d$-dimensional polytope with $n$ facets
that can lie on a \emph{monotone path}, i.e., on a path along edges
that is strictly increasing with respect to a linear objective
function.

McMullen's~1971 \emph{Upper Bound Theorem}~\cite{McMullen71} (claimed
by Motzkin~\cite{Motzkin57} in~1957) states that the maximal number
$M_{\rm ubt}(d,n)$ of vertices that any $d$-dimensional polytope with
$n$~facets can have is achieved by the polars $C_d(n)^\Delta$ of
cyclic $d$-polytopes with $n$~facets.

The Upper Bound Theorem yields, for all $n>d\ge2$, the inequality
\begin{equation}\label{eq:ubt}
  M(d,n)\ \le\ M_{\rm ubt}(d,n),
\end{equation}
but from this it is not clear whether equality always holds, that is,
if for all $n>d\ge2$ one can construct a simple polar-to-neighborly
$d$-polytope with $n$~facets that admits a monotone Hamilton path with
respect to a linear objective function. Equality in~\eqref{eq:ubt}
\emph{is} known in the cases $d\leq 3$~and~$n\leq d+2$.

However, in~\cite{PfeifleZiegler04a} we show that in fact
$M(6,9)<M_{\rm ubt}(6,9)$: there exists \emph{no} realization of the
(combinatorially unique) polar-to-neighborly $6$-polytope
$C_6(9)^\Delta$ with $9$~facets and $30$~vertices that admits such a
monotone Hamilton path.

For the parameters $d=4$, $n=8$, one can show using the same
(basically combinatorial) methods that there is also no realization
of~$C_4(8)^\Delta$ with a monotone Hamilton path --- but as we will
show here, there are other dual-to-neighborly but not dual-to-cyclic
$4$-dimensional polytopes with $8$~facets that admit a realization
with a monotone path through all vertices.

In fact, in this paper we prove considerably more: we provide a
geometric construction that shows that the inequality~\eqref{eq:ubt}
is tight in dimension~$d=4$ for all~$n\ge5$.

\begin{maintheorem} \label{thm:main} 
For each integer $m\geq 0$, there exists a simple polar-to-neighborly
$4$-dimensional polytope $Q_m$ with $n=m+5$ facets and a linear
objective function $f:\R^4\to\R$, such that the orientation induced
by~$f$ on the 1-skeleton of~$Q_m$ admits a monotone Hamilton
path. Therefore, 
\[
    M(4,n) \ = \ M_{\rm ubt}(4,n) \ = \ \tfrac12\,n(n-3).
\]
In other words, the maximal number $M(4,n)$ of vertices on a strictly
monotone path in the graph of a $4$-dimensional polytope with $n$~facets
equals the maximal number of vertices that such a polytope can have
according to the Upper Bound Theorem.
\end{maintheorem}

An interesting feature used in our proof is that for $m\ge3$, the
(polar-to-)neighborly polytopes~$Q_m$ are \emph{not} polar to cyclic
ones. In fact, exhaustive enumeration shows that already the graph
of~$C_4(8)^\Delta$ does not satisfy a combinatorial condition
necessary for the existence of an monotone path, namely, it does not
admit a Hamilton AOF Holt-Klee orientation~\cite{Holt-Klee99}.  This
is also true for the graphs of the polytopes $C_4(n)^\Delta$ for $8\le
n\le 12$; we conjecture that the graphs of $C_4(n)^\Delta$ for
all~$n\ge8$ admit no Hamilton AOF Holt-Klee orientation.

\smallskip The structure of the paper is as follows: We first give an
explicit description, reminiscent of Gale's Evenness Criterion for
polar-to-cyclic polytopes, of the combinatorial structure of a
family $\{Q^d_m:d\ge4\text{ even},\, m\ge 0\}$ of simple
(polar-to-)neighborly $d$-dimensional polytopes with $m+d+1$~facets 
(Sections~\ref{sec:results}~and~\ref{sec:comb}).  For $d=4$, we then use
this description to specify a Hamilton path $\pi_m$ on
each~$Q_m:=Q^4_m$ (Section~\ref{sec:ham-path}).  In
Section~\ref{sec:realizing-paths}, we start with a monotone
path~$\pi_0$ on a certain realization of the $4$-simplex~$Q_0$, and
for~$m\ge0$ inductively realize the polytope~$Q_{m+1}$ in such a way
that the path~$\pi_{m+1}$ is strictly monotone with respect to a
suitable objective function (Theorem~\ref{thm:detailed}). We proceed
in three steps: First, we position~$Q_m$ in a suitable way with
respect to the standard coordinates on~$\R^4$
(Section~\ref{subsec:step1}). We then find a ``cutting
plane''~$H_{m+1}$ such that the polytope $Q_m\cap H_{m+1}^{\ge0}$ has
the right combinatorial type (Section~\ref{subsec:step2}). Finally, we
complete the construction in Section~\ref{subsec:step3} by applying a
projective transformation~$\psi$ to~$\R^4$ such that the path
$\psi(\pi_m)$ on $Q_{m+1}:=\psi(Q_m\cap H_{m+1}\pos)$ is strictly
monotone with respect to the objective function $f:\R^4\to\R$,
$\x\mapsto x_4$.

\section{Main results}\label{sec:results}

\begin{theorem}[modified Gale's Evenness Criterion] \label{thm:mod-gale}
For each $m\ge0$ and even $d\ge4$, the  following sets correspond to
the vertices of a combinatorial type $\Qt^d_m$ of a simple $d$-dimensional
polar-to-neighborly polytope with $n=m+d+1$ facets. 

\smallskip
  \begin{compactitem}[\quad$\triangleright$]    
    \item\emph{Type 1.} The union of one ``triplet with a hole'' and
      $d/2-1$ pairs of indices
      \[ 
      \{j_1,\,j_1+2\}\;\cup\;\{j_2,\,j_2+1\}\;\cup\;\cdots
      \;\cup\;\{j_{d/2},j_{d/2}+1\},
      \] 
      where $1\le j_1<n-d+1$, $j_1+3\le j_2$, $j_k+2 \le j_{k+1}$ for
      $2\le k \le d/2-1$, and $j_{d/2}<n$.
    \item\emph{Type 2a.} The union of one triplet, the
      singleton~$\{n\}$, and $d/2-2$ pairs of indices
      \[
      \{j_1,\,j_1+1,\,j_1+2\}\;\cup\;\{j_2,\,j_2+1\}\;\cup\; \cdots
      \;\cup\;\{j_{d/2-1},j_{d/2-1}+1\}\;\cup\;\{n\},
      \]
      where $1\le j_1<n-d+1$, $j_1+3\le j_2$, $j_k+2 \le j_{k+1}$ for
      $2\le k \le d/2-2$, and $j_{d/2-1}<n-1$.
    \item\emph{Type 2b.} The union of $d/2$ pairs of indices
      \[ 
      \{1,2\}\;\cup\;\{j_1,\,j_1+1\}\;\cup\;\cdots 
      \;\cup\;\{j_{d/2-1},j_{d/2-1}+1\},
      \] 
      where $3\le j_1$, $j_k+2 \le j_{k+1}$ for $2\le k \le d/2-2$,
      and~$j_{d/2-1}<n$.
  \end{compactitem}
\end{theorem}

  \begin{figure}[h]
    \centering
    \vspace*{-4cm}
        \psfrag{j1}{\small$j_1$} \psfrag{j2}{\small$j_2$}
        \psfrag{jd2}{\small$j_{d/2}$} \psfrag{jf-1}{\small$j_{d/2-1}+1$}
        \psfrag{1}{\small1}\psfrag{n}{\small$n$}
        \psfrag{Type 1}{Type 1}
        \psfrag{Type 2a}{Type 2a}
        \psfrag{Type 2b}{Type 2b}
        \includegraphics[width=.7\linewidth,angle=-90]{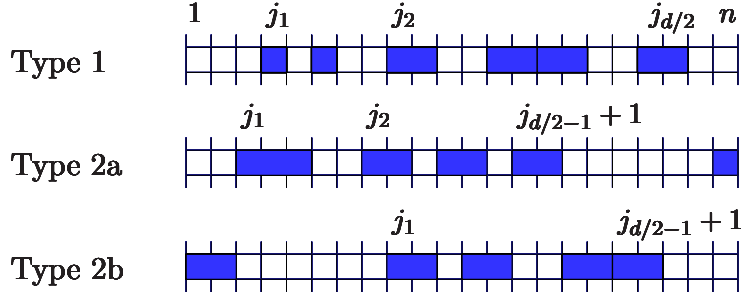} 
        \vspace*{-4cm}
    \caption{The vertex-facet incidences of the polytopes $\Qt^d_m$ are
      obtained from these patterns by fixing the dark boxes, and
      sliding the lighter boxes between $1$~and~$n$ without overlap.
      For Type~1, the box $\{i,i+2\}$ must be regarded as one rigid
      unit.}
    \label{fig:vif-Q4}
  \end{figure}

\begin{remark} \label{rem:ptn}
  If we accept for the moment the existence of the polytopes $\Qt^d_m$,
  it is easy to verify that they are polar-to-neighborly by counting
  the number of vertices using Figure~\ref{fig:vif-Q4}:
  \begin{eqnarray*}
    f_0(\Qt^d_{n-d-1}) &=& 
          \underbrace{ {n-2-(d/2-1) \choose d/2} }_{\text{Type 1}} +
          \underbrace{ {n-2-(d/2-2)-1 \choose d/2-1 } }_{\text{Type 2a}}
          +\underbrace{ {n-2-(d/2-1) \choose d/2-1 } }_{\text{Type 2b}} \\
              &=& {n-1-d/2 \choose d/2} + 2{n-1-d/2\choose d/2-1} \\[1ex]
              &=& {n-d/2 \choose d/2} + {n-1-d/2\choose d/2-1}, 
  \end{eqnarray*}
  which is the number of vertices of a simple polar-to-neighborly
  $d$-polytope with $n=m+d+1$ facets, since $d$ is assumed even. By
  \cite[Chapter 8]{Ziegler98}, \emph{any} polytope with that many
  vertices is polar-to-neighborly.~$\qed$
\end{remark}

From now on, we will always write $\Qt_m := \Qt^4_m$.

\begin{figure}[htbp]
  \centering \includegraphics[width=.64\linewidth]{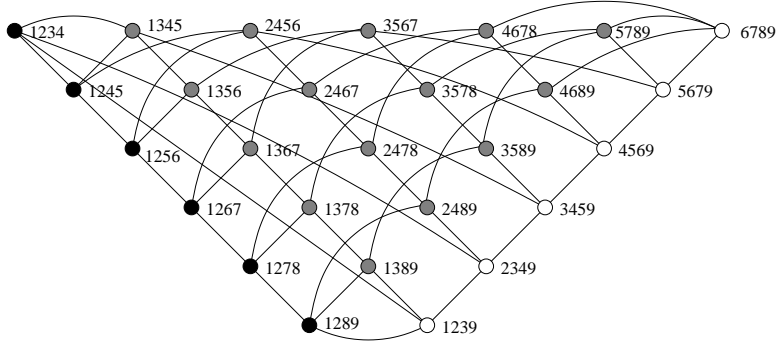}
  \caption{Graph of the $4$-polytope $\Qt_4$ with $n=9$
    facets. Vertices of type 1, 2a, and 2b are drawn in gray, white,
    and black, respectively. Each vertex is labelled with the facets
    it is incident~to.}
  \label{fig:q49-graph}
\end{figure}

\begin{proposition} \label{prop:path}
  Each polytope $\Qt_m$ admits a Hamilton path $\pit_m$ in its
  graph that induces an AOF-orientation
  (cf.~Figure~\ref{fig:q49-graph+path} and Definition~\ref{def:aof} below).
\end{proposition}

\begin{figure}[htbp]
  \centering
  \psfrag{T0}{$T^0_4=\{\beta_4\}$} \psfrag{T1}{$T^1_4=\{\alpha_4\}$}
  \psfrag{T2}{$T^2_4$} \psfrag{T3}{$T^3_4$}
  \psfrag{T4}{$T^4_4$} 
  \psfrag{p}{} \psfrag{r}{$\tau_4$} \psfrag{z}{$\omega_4$}
  \includegraphics[width=.6\linewidth]{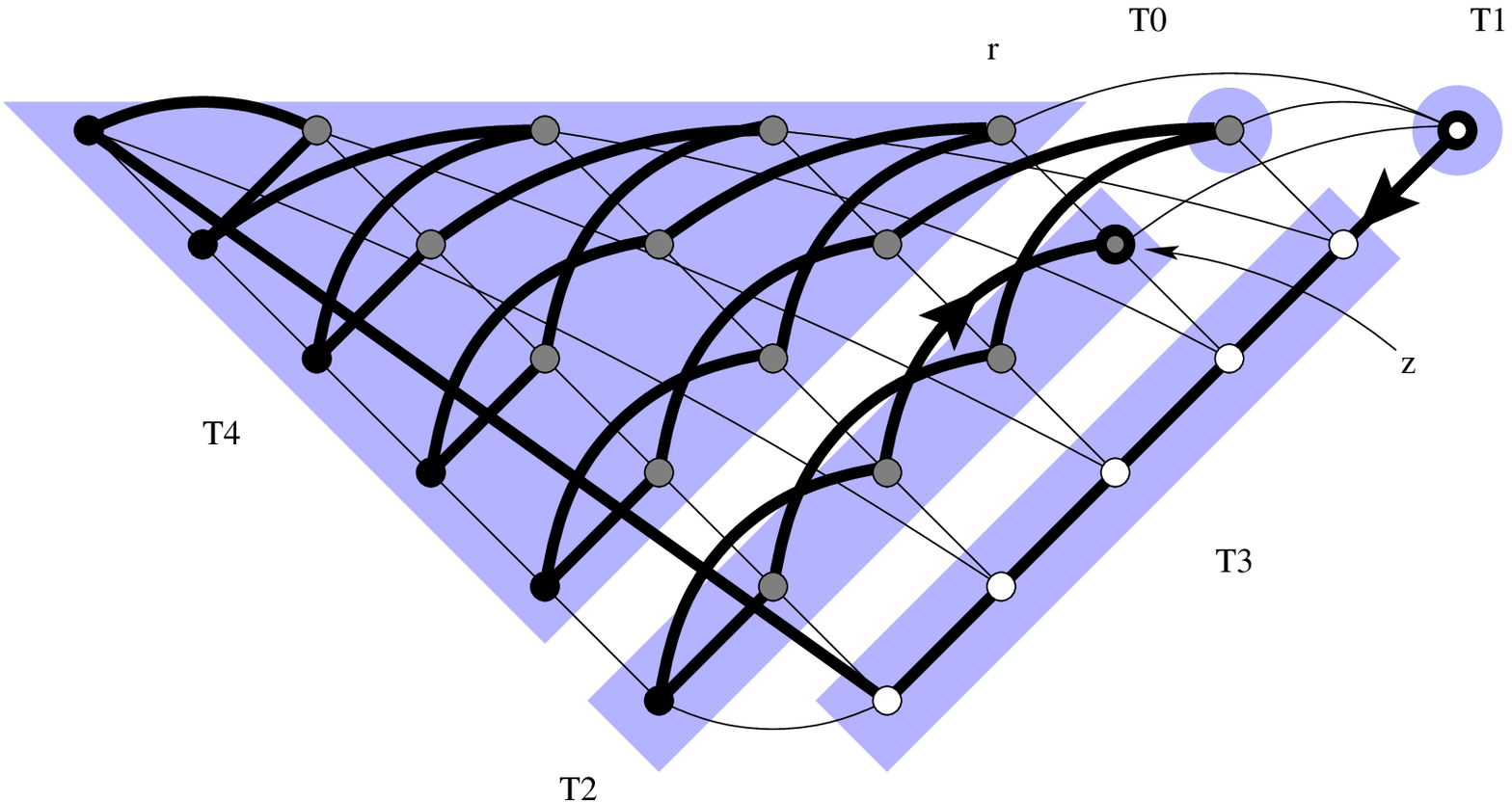}\qquad
  \includegraphics[width=.33\linewidth]{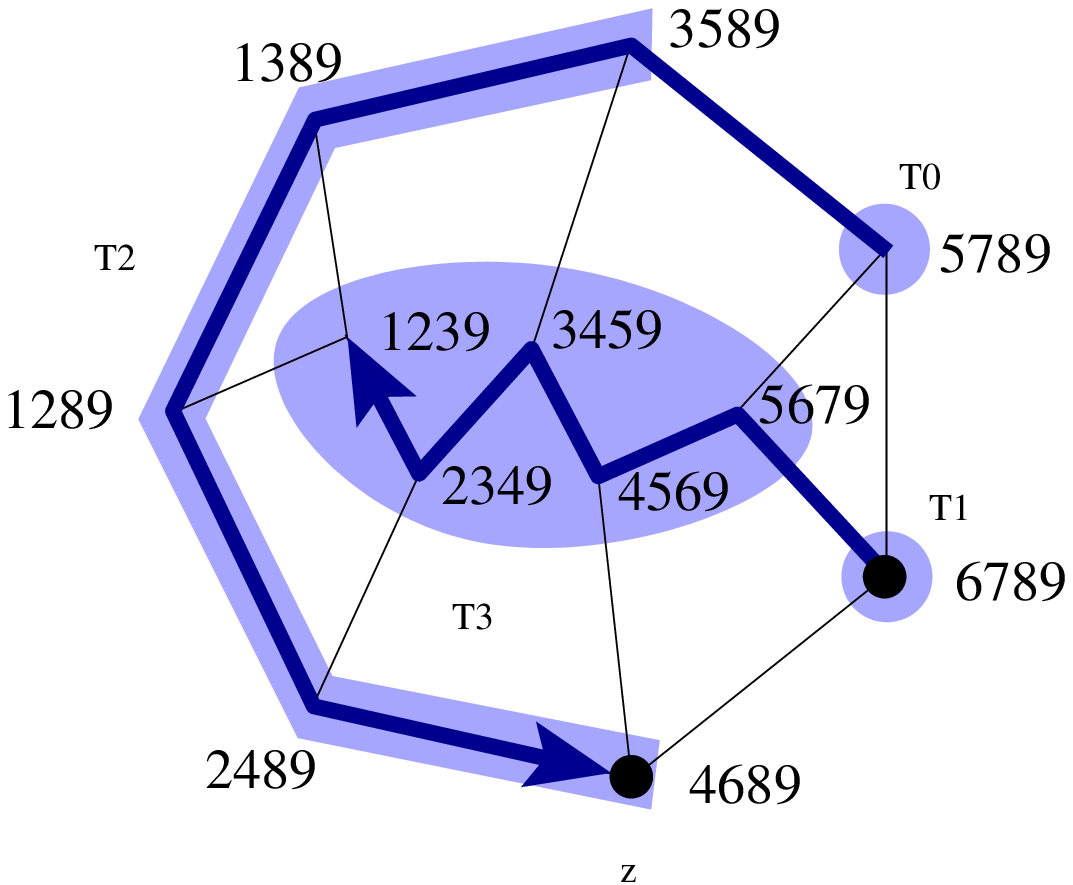}
  \caption{\emph{Left:} Graph of $\Qt^4_4$. The partition of the
    vertices into the tips $T^0$, $T^1$, \dots, $T^4$ is shown, along
    with the Hamilton path $\pit_4$ (bold). The source~$\alpha_m$
    is labeled $\{n-3,n-2,n-1,n\}$, and the sink
    $\omega_m=\{n-5,n-3,n-1,n\}$.  See Convention~\ref{conv:vertices}
    for the labels of the other marked vertices.  \emph{Right:} The
    facet~$F^3_4$ with the restriction of~$\pit_4$ to it.}
  \label{fig:q49-graph+path}
\end{figure}

\begin{remark}\label{rem:order-pi}
  The crucial property for our realization construction is that the
  path~$\pit_m$ begins in a certain facet~$F_m^3$ of the
  polytope~$Q_m$ (defined below), traverses the rest of~$Q_m$, and
  then returns to~$F_m^3$ (cf.~Figure~\ref{fig:q49-graph+path}). This
  permits us to add new vertices to the beginning and end of~$\pit_m$
  by modifying only the facet~$F_m^3$.
\end{remark}

\begin{theorem} \label{thm:detailed} 
  There exists a family $\{Q_m:m\ge0\}$ of special realizations of the
  combinatorial types~$\Qt_m$, in which each Hamilton path
  $\pi_m$~visits the vertices of~$Q_m$ in the order given by
  increasing $x_4$-coordinate.  This family may be realized
  inductively starting from the $4$-simplex~$Q_0$ in such a way that
  for all~$m\ge0$, a realization of $Q_{m+1}$ with a monotone
  Hamilton path~$\pi_{m+1}$ may be obtained from \emph{any}
  realization of~$Q_m$ with such a path~$\pi_m$.
\end{theorem}

\section{Constructing the combinatorial types $\Qt^d_m$}
\label{sec:comb}

\subsection{Facet splitting}

We will prove Theorem \ref{thm:mod-gale} using Barnette's technique of
\emph{facet splitting}~\cite{Barnette81}. Put briefly, for each even
$d\ge4$ we will inductively construct a family
$\{(\Qt^d_m,\F_m):m\ge0\}$, where each $\Qt^d_m$ is the combinatorial
type of a simple $d$-dimensional polytope with $m+d+1$ facets, and
$\F_m$ is a flag of faces on~$\Qt^d_m$ (to be defined shortly). We then
use~$\F_m$ to find a ``good'' oriented hyperplane~$H_{m+1}$ in general
position with respect to the vertices of~$\Qt^d_m$, and set
$\Qt^d_{m+1}:=\Qt^d_m\cap H_{m+1}\pos$.

\begin{definition}
  Let $P$ be a $d$-dimensional simple polytope. A \emph{flag of faces}
  on $P$ is a chain
\begin{equation}\label{eq:flag}
    \F:\quad \emptyset=F^{-1}\ \subset\ F^0\ \subset\
    F^1\ \subset\ \cdots\ \subset\  F^d=P
\end{equation}
  of faces of $P$ such that $\dim F^i=i$ for~$i=0,1,\dots,d$.  The
  \emph{$i$-th tip} of a flag $\F$ is $T^i:=\vertices F^i\setminus
  \vertices F^{i-1}$, for $0\le i\le d$. We say that the tip~$T^i$ is
  \emph{even} resp.~\emph{odd} according to the parity of~$i$.
  Moreover, for $0\le k\le d$ we set
  \[
      T_{\rm even}^{\le k} \ = \ 
        \bigcup_{0\le e \le k\atop e\text{ even}} T^e
      \qquad\text{and}\qquad
      T_{\rm odd}^{\le k}\ = \ 
        \bigcup_{1 \le o \le k\atop o\text{ odd}} T^o.
  \]
\end{definition}

\begin{lemma}\label{lem:splittingplane}
  Let $P$ be a simple $d$-dimensional polytope with $n$ facets, and
  $\F$ a flag of faces as in~\eqref{eq:flag}. Then there exists an
  affine oriented hyperplane $H$ in general position with respect
  to~$P$ such that $T_{\rm even}^{\le d}\subset H^+$ and $T_{\rm
  odd}^{\le d}\subset H^-$.  In particular, $P\cap H\pos$ is a simple
  $d$-polytope with $n+1$ facets.
\end{lemma}

\begin{proof}  Pick an oriented point $\{v\}=H^0\subset\relint
F^1$ such that $T^0\in (H^0)^+$. Inductively, for $1\le k\le d-1$,
if we have already chosen an oriented $(k-1)$-dimensional
affine subspace~$H^{k-1}$ in~$\aff{F^k}$ such that
\begin{equation} \label{eq:star}
    T_{\rm even}^{\le k} \ \subset\ (H^{k-1})^+
    \qquad\text{and}\qquad
    T_{\rm odd}^{\le k} \ \subset\ (H^{k-1})^-, 
\end{equation}
we take a $k$-plane~$H^k$ that initially coincides with~$\aff{F^k}$,
and orient it in such a way that $T^{k+1}$ lies in~$(H^k)^+$ if $k+1$
is even, respectively in~$(H^k)^-$ if $k+1$ is odd. Now we
rotate~$H^k$ by a sufficiently small amount around~$H^{k-1}$ in such a
way that $T_{\rm even}^{\le k}\subset (H^k)^+$.  Then
\eqref{eq:star}~even holds with $k$~replaced by~$k+1$. By
construction, the hyperplane $H:=H^{d-1}$ is in general position with
respect to~$P$.
\end{proof}

\begin{definition} \label{def:Qt}
  The family $\{(\Qt^d_m, \F_m):m\ge0\}$ of $d$-dimensional polytopes
  $\Qt^d_m$ equipped with flags $\F_m$ of faces is defined in the
  following way:
\begin{compactenum}[(a)]
\item $\Qt^d_0$ is the combinatorial type of the $d$-simplex
$\conv\{v_1,v_2,\dots, v_{d+1}\}$. The flag $\F_0$ is defined by setting
$F_0^i := \conv\big\{v_1,\,v_2,\dots,v_{i+1}\big\}$ for
$i=0,1,\dots,d$.  Then
\begin{equation}\label{eq:toi}
  T_0^i \ :=\
  \vertices\big(F_0^i\big)\setminus\vertices\big(F_0^{i-1}\big)
  \ =\ \{v_{i+1}\} \qquad \text{for }\ i=0,1,\dots,d.
\end{equation}

\item For $m\ge0$, let $H=H_{m+1}$ be the oriented hyperplane given by
applying Lemma \ref{lem:splittingplane} to $P=\Qt^d_m$ and $\F=\F_m$,
and set $\Qt^d_{m+1} := \Qt^d_m\cap H_{m+1}\pos$ and (cf.\ 
Figure \ref{fig:newtips})
 \begin{eqnarray*}\label{eq:newtips}
      T_{m+1}^0 &:=& \vertices \big(\conv(T_m^1 \cup T_m^2) \cap H_{m+1}\big),\\
      T_{m+1}^1 &:=& \vertices \big(\conv(T_m^0 \cup T_m^1) \cap H_{m+1}\big),\\
      T_{m+1}^j &:=& \vertices \Big(\conv\bigg(T_m^{j+1} \cup
      \displaystyle\bigcup_{0\le k<j \atop k+j=0\bmod 2}
      T_m^k \bigg) \cap H_{m+1} \Big)  \quad\text{for }\; j=2,3,\dots,d-1,\\
      T_{m+1}^d &:=& \bigcup_{0\le k\le d/2} T_m^{2k}.
    \end{eqnarray*}

\begin{figure}[htbp]
  \centering
  \input{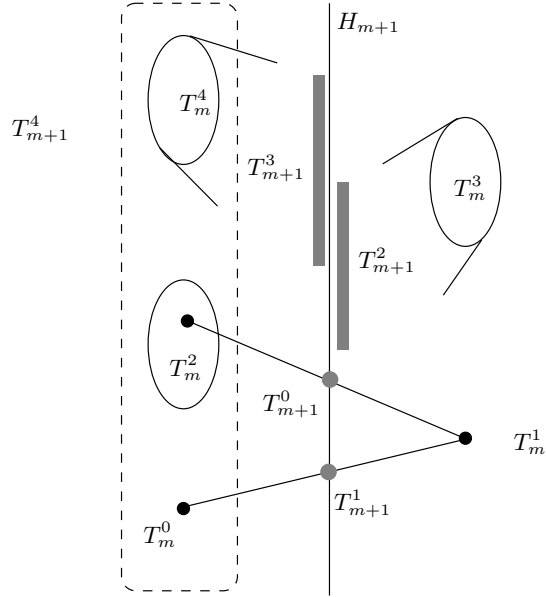}
  \caption{New tips in the case $d=4$.}
  \label{fig:newtips}
\end{figure}

\noindent The flag $\F_{m+1}$ is now defined by $F_{m+1}^j:=\bigcup_{i=0}^j
T^i_{m+1}$ for $j=0,1,\dots,d$. Moreover, put
  \[
      T_{\rm even}^{\le k}(m) \ = \ 
        \bigcup_{0\le e \le k\atop e\text{ even}} T^e_m
      \qquad\text{and}\qquad
      T_{\rm odd}^{\le k}(m)\ = \ 
        \bigcup_{1 \le o \le k\atop o\text{ odd}} T^o_m.
  \]
\end{compactenum}
\end{definition}

\begin{remark} \
  \begin{compactenum}[(a)] 
  \item The polytopes $C_d(n)^\Delta$ arise by
    exchanging the definitions of~$T^0_{m+1}$ and~$T^1_{m+1}$.
    
  \item  All new vertices arise
    as the intersection of~$H_{m+1}$ with some edge $\conv\{v,w\}$
    of~\smash{$\Qt^d_m$}, where $v$~and~$w$ lie in tips of different
    parity. Furthermore, all vertices of~\smash{$\Qt^d_m$} belonging to even
    tips are also vertices of~\smash{$\Qt^d_{m+1}$}, and vertices in odd tips
    disappear.
  \end{compactenum}
\end{remark}

\begin{proposition}\label{prop:q_m:vertices}
  For each $m\ge0$, the following is true for the pair~$(\Qt^d_m, \F_m)$:
  \begin{compactenum}[(a)] 
  \item For all $i,j\in\N$ with $0\le i<j\le d$ and $i+j=1\bmod 2$ and
    all $v\in T^i_m$, there is exactly one~$w\in T^j_m$ such that
    $\conv\{v,w\}\in\sk^1(\smash{\Qt^d_m})$. This gives rise to bijections
    $T^{\le k}_{\rm even}(m)\cong T^k_{m+1}$ for odd $0<k<d$
  resp. $T^{\le k}_{\rm odd}(m)\cong T^k_{m+1}$ for even $0\le k\le d$.
  \item $|T_m^e|=|T_m^{e+1}|={e/2+m\choose m}$ for even $e=0,2,\dots, d-2$,
    and $|T_m^d|={d/2+m\choose m}$. This proves again that 
    \smash{$\Qt^d_m$}~is polar-to-neighborly.
  \end{compactenum}
\end{proposition}

\begin{proof} (a) This follows because $v$~lies in
$F_m^{j-1}=\bigcup_{i=0}^{j-1} T^i_m$, and $\conv(F_m^{j-1})$~is a
$(j-1)$-dimensional face of the simple polytope
$\conv(F_m^j)=\conv(F_m^{j-1}\cup T_m^j)$.

\noindent (b)  We proceed by induction, and can assume that
the assertion holds for $m\ge0$. From the bijections in part (a), we
conclude for all even $e=0,2,\dots, d-2$ that
\[
   |T_{m+1}^{e+1}| \;=\; |T_{m+1}^e|\;=\;\sum_{i=0\atop i\text{ even}}^e |T_m^i|\;=\;
   \sum_{k=0}^{e/2} |T_m^{2k}| \;=\; \sum_{k=0}^{e/2} {k+m\choose m}
   \;=\; {e/2+m+1\choose m+1}.
\]
The calculation for $|T^d_{m+1}|$ is similar. The fact that
\smash{$\Qt^d_m$}~is polar-to-neighborly follows by the same argument as
in Remark~\ref{rem:ptn}, since
\begin{eqnarray*}
  f_0(\Qt^d_m) &=& \sum_{k=0}^{d/2}{k+m\choose m} +
  \sum_{k=0}^{\lfloor(d-1)/2\rfloor}{k+m\choose m} \\
  &=& {m + d/2 + 1 \choose d/2} + 
      {m + \lfloor(d-1)/2\rfloor + 1 \choose
        \lfloor(d-1)/2\rfloor}\\[1ex]
  &=& {n - d/2 \choose d/2} + {n - \lceil(d-1)/2\rceil - 1 \choose
    \lfloor(d-1)/2\rfloor}.
\end{eqnarray*}

\end{proof}

\subsection{Combinatorics of the family \boldmath $\Qt^d_m$}

\begin{convention} \label{conv:labels}
We introduce labelings to make the combinatorics of the $\Qt_m$
explicit:

  \begin{compactenum}[(a)] 
  \item For any labeling of the
        facets of a simple $d$-polytope~$P$ with labels
        in~$[n]:=\{1,2,\dots,n\}$,  let $\lambda:\vertices P\to
        {[n]\choose d}$ assign to each vertex~$v$ of~$P$ the set of
        labels of all facets that $v$~is incident to. We identify a
        vertex~$v$ with its label~$\lambda(v)$.
    
  \item The facets of the $d$-simplex~$\Qt^d_0$ on the vertex set
    $\{v_1,v_2,\dots,v_{d+1}\}$ are labeled in such a way that
    $v_1\equiv\lambda(v_1)=[d+1]\setminus\{2\}$,
    $v_2=[d+1]\setminus\{1\}$, and $v_j=[d+1]\setminus\{j\}$
    for~$j=3,4,\dots,d+1$ (cf.\ Figure~\ref{fig:4-simplex}).
    
  \item The ``new'' facet~$\Qt^d_m\cap H_{m+1}$ of $\Qt^d_{m+1}$ is
    labelled $m+d+2$.  

  \end{compactenum}
\end{convention}

\begin{figure}[htbp]
  \centering
  \renewcommand{\arraystretch}{-2}
  \begin{tabular}[t]{rl}
    \begin{tabular}[c]{r}
      $\{v_5\} \;=\;T_0^4$
    \end{tabular}
    \begin{tabular}[c]{>{\scriptsize}r@{\;\;}l}
      2b & 12|34\\
    \end{tabular}
    \\
    &
    \begin{tabular}[c]{r@{\;\;}>{\scriptsize}l}
      123|5 & 2a
    \end{tabular}
    \begin{tabular}[c]{l}
      $T_0^3 \;=\; \{v_4\}$
    \end{tabular}     
    \\
    \begin{tabular}[c]{r}
      $\{v_3\} \;=\; T_0^2$
    \end{tabular}
    \begin{tabular}[c]{>{\scriptsize}r@{\;\;}l}
      2b & 12|45
    \end{tabular}
    \\
    &
    \begin{tabular}[c]{r@{\;\;}>{\scriptsize}l}
      234|5 & 2a
    \end{tabular}
    \begin{tabular}[c]{l}
      $T_0^1 \;=\; \{v_2\}$
    \end{tabular}     
    \\
    \begin{tabular}[c]{r}
      $\{v_1\} \;=\; T_0^0$
    \end{tabular}
    \begin{tabular}[c]{>{\scriptsize}r@{\;\;}l}
      \phantom{2}1  & 13|45
    \end{tabular}
  \end{tabular}
  
  \caption{The labeling of the vertices of the $4$-simplex~$\Qt_0$ according to
    Convention~\ref{conv:labels}(b). Also shown is the classification
    of the vertices into types 1, 2a, 2b as in
    Proposition~\ref{prop:vif-qd}. }
  \label{fig:4-simplex}
\end{figure}

\begin{proposition}\label{prop:vif-qd} 
  Let $m\ge0$ and $n=m+d+1$.  A vertex $v$ of
  $\Qt^d_m$ lies in $T^i_m$ exactly if 
  \[
        {\textstyle\max_n\vb} \;:=\; \max\big([n]\setminus v\big) \;=\;
        \begin{cases}
          m+2   & \text{for } i=0,\\ 
          m+1   & \text{for } i=1,\\
          m+i+1 & \text{for } 2\le i \le d.
        \end{cases}
  \]
\end{proposition}

\begin{proof}  This is true for $m=0$ by~\eqref{eq:toi} and
  Convention~\ref{conv:labels}, see also
  Figures~\ref{fig:4-simplex}~and~\ref{fig:dim4}.  For $m>0$ and~$2\le
  i \le d-1$, the statement follows because any vertex~$\vt\in
  T_m^{i}$ is of the form $\vt=\conv\{v,w\}\cap H_m\equiv(v\cap
  w)\cup\{n\}$ for some $v\in T_{m-1}^k$ and $w\in T_{m-1}^{i+1}$ with
  $k\le i$. But then by induction,
\[
   {\textstyle \max_{n-1}\vb \;<\; \max_{n-1}\wb} \;=\; 
   (m-1) + (i+1) + 1 \;=\; m+i+1,
\]
so $\max\big([n]\setminus\vt\big)=m+i+1$ as required. The case~$i=d$
follows directly from Definition~\ref{def:Qt}, and the cases $i=0,1$
are checked similarly.
\end{proof}

\begin{proof}[Proof of Theorem \ref{thm:mod-gale}]
The existence of the family $\{(\Qt^d_m, \F_m):m\ge0\}$ follows from
Lemma~\ref{lem:splittingplane}. Using
Propositions~\ref{prop:q_m:vertices}~and~\ref{prop:vif-qd}, it is
somewhat tedious but elementary to verify that for all $m\ge0$, the
vertices of $\Qt^d_m$ are of the given types. More precisely, all
vertices of~$T_{\rm even}^{\le d}(m)$ are of type~1 or~2b, and $T_{\rm
odd}^{\le d-1}(m)$ is made up entirely of vertices of type~2a,
cf.~Figure~\ref{fig:dim4}.
\end{proof}


\begin{figure}[htbp]
  \renewcommand{\arraystretch}{1}
  \centering\small
  \begin{minipage}[b]{.49\linewidth}
    \vspace{0pt}
    \centering
    \begin{tabular}[t]{c|c}
      \begin{tabular}[c]{r}
        $T_1^4$
      \end{tabular}
      \begin{tabular}[c]{>{\scriptsize}r@{\quad}l}
        2b & 12|34\\
        2b & 12|45\\
        1  & 13|45
      \end{tabular}
      \hspace{-1ex}
      \begin{tabular}[c]{r}
        \scriptsize 6
      \end{tabular}
      \\
      &
      \begin{tabular}[c]{r}
        \scriptsize 5
      \end{tabular}
      \hspace{-1ex}
      \begin{tabular}[c]{r@{\quad}>{\scriptsize}l}
        123|6 & 2a\\
        234|6 & 2a
      \end{tabular}
       \begin{tabular}[c]{r}
        $T_1^3$
      \end{tabular}     
      \\
      \begin{tabular}[c]{r}
        $T_1^2$
      \end{tabular}
      \begin{tabular}[c]{>{\scriptsize}r@{\quad}l}
        2b & 12|56\\
        1  & 13|56
      \end{tabular}
      \hspace{-1ex}
      \begin{tabular}[c]{r}
        \scriptsize 4
      \end{tabular}
      \\
      &
      \begin{tabular}[c]{r}
        \scriptsize 2
      \end{tabular}
      \hspace{-1ex}
      \begin{tabular}[c]{r@{\quad}>{\scriptsize}l}
        345|6 & 2a
      \end{tabular}
       \begin{tabular}[c]{r}
        $T_1^1$
      \end{tabular}     
      \\
       \begin{tabular}[c]{r}
        $T_1^0$
      \end{tabular}
      \begin{tabular}[c]{>{\scriptsize}r@{\quad}l}
        \ 1  & 24|56
      \end{tabular}
      \hspace{-1ex}
      \begin{tabular}[c]{r}
        \scriptsize 3
      \end{tabular}
    \end{tabular}
    \mbox{}
  \end{minipage} 
  \hfill
  \begin{minipage}[b]{.49\linewidth}
    \vspace{0pt}
    \centering
    \begin{tabular}[t]{c|c}
      \begin{tabular}[c]{r}
        $T_2^4$
      \end{tabular}
      \begin{tabular}[c]{>{\scriptsize}r@{\quad}l}
        2b & 12|34\\
        2b & 12|45\\
        1  & 13|45\\
        2b & 12|56\\
        1  & 13|56\\
        1  & 24|56
      \end{tabular}
      \hspace{-1ex}
      \begin{tabular}[c]{r}
        \scriptsize 7
      \end{tabular}
      \\
      &
      \begin{tabular}[c]{r}
        \scriptsize 6
      \end{tabular}
      \hspace{-1ex}
      \begin{tabular}[c]{r@{\quad}>{\scriptsize}l}
        123|7 & 2a\\
        234|7 & 2a\\
        345|7 & 2a
      \end{tabular}
       \begin{tabular}[c]{r}
        $T_2^3$
      \end{tabular}     
      \\
      \begin{tabular}[c]{r}
        $T_2^2$
      \end{tabular}
      \begin{tabular}[c]{>{\scriptsize}r@{\quad}l}
        2b & 12|67\\
        1  & 13|67\\
        1  & 24|67
      \end{tabular}
      \hspace{-1ex}
      \begin{tabular}[c]{r}
        \scriptsize 5
      \end{tabular}
      \\
      &
      \begin{tabular}[c]{r}
        \scriptsize 3
      \end{tabular}
      \hspace{-1ex}
      \begin{tabular}[c]{r@{\quad}>{\scriptsize}l}
        456|7 & 2a
      \end{tabular}
       \begin{tabular}[c]{r}
        $T_2^1$
      \end{tabular}     
      \\
       \begin{tabular}[c]{r}
        $T_2^0$
      \end{tabular}
      \begin{tabular}[c]{>{\scriptsize}r@{\quad}l}
        \ 1 & 35|67
      \end{tabular}
      \hspace{-1ex}
      \begin{tabular}[c]{r}
        \scriptsize 4
      \end{tabular}
   \end{tabular}
   \mbox{}
  \end{minipage}

  $\underbrace{\hspace{.45\linewidth}}_{\displaystyle\Qt_1}$\qquad
  $\underbrace{\hspace{.45\linewidth}}_{\displaystyle\Qt_2}$
  \caption{Vertex labels in the polytopes $\Qt_1$ (left) and
    $\Qt_2$ (right). Also shown are the type (outside) of each vertex~$v$ and
    the value of $\max_n\bar v$ (inside).}
  \label{fig:dim4}
\end{figure}

\section{A Hamilton path $\pit_m$ 
        that induces an AOF-orientation on $\Qt_m$ }
\label{sec:ham-path}

\begin{definition} \label{def:aof}
Let $P$~be a simple $d$-polytope. An acyclic orientation of the graph
of~$P$ that has a unique sink in each face (including $P$~itself) is
called an \emph{AOF-orientation}
on~$P$. For any orientation~$\calO$ of the graph of~$P$ and $0\le k\le
d$, denote by $h_k(\calO)$ the number of vertices of in-degree~$k$
in~$\calO$.

\end{definition}

\begin{proposition} \label{prop:hk} (see
  e.g.~\cite[Chap.~8.3]{Ziegler98} and \cite{Joswig-Kaibel-Koerner02})
  An acyclic orientation~$\calO$ of the graph of a simple
  $d$-polytope~$P$ is an AOF-orientation if and only if the $h$-vector
  of~$P$ coincides with the vector $\big(h_0(\calO),h_1(\calO),\dots,
  h_d(\calO)\big)$. \hfill$\Box$
\end{proposition}

\begin{proof}[Proof of Proposition \ref{prop:path}]
  By inspection of Figures~\ref{fig:q49-graph}
  and~\ref{fig:q49-graph+path}, the algorithm of Figure~\ref{fig:path}
  yields a Hamilton path $\pit_m$ in the graph of~$\Qt_m$. Note that
  $\pit_m$ passes through $T_m^1$, $T_m^3$, $T_m^4$, $T_m^0$,
  and~$T_m^2$, in this order (cf.~Remark~\ref{rem:order-pi}).

\begin{figure}[htbp]
  \begin{compactenum}
  \item \emph{``Odd stage''.}
    \begin{minipage}[t]{.8\linewidth}
      \begin{tabbing}
        \qquad\=\qquad\=\qquad\=\hspace{5cm}\=\kill
        \> \textbf{for} $i$ \textbf{from} $n-3$ \textbf{to} $1$ \textbf{do} \\
        \>\> visit $\{i,\,i+1,\,i+2,\,n\}$;
      \end{tabbing}
    \end{minipage}
  \item \emph{``Even stage''.}
    \begin{minipage}[t]{.8\linewidth}
      \begin{tabbing}
        \qquad\=\qquad\=\qquad\=\hspace{5cm}\=\kill
        \> \textbf{for} $j$ \textbf{from} $3$ \textbf{to} $n-1$ \textbf{do} \\
        \>\> $i:=j-3$;\\
        \>\> \textbf{while} $i\ge1$ \textbf{do}  \>\> \emph{``down'' phase} \\
        \>\>\> visit $\{i,\,i+2,\,j,\,j+1\}$;\\
        \>\>\> $i:=i-2$;\\
        \>\> visit $\{1,\,2,\,j,\,j+1\}$;\\
        \>\> \textbf{if} $j$ is even \textbf{then} $i:=2$; \textbf{else} $i:=1$;\\
        \>\> \textbf{while} $i\le j-4$ \textbf{do}\>\> \emph{``up''
          phase}\\
        \>\>\> visit $\{i,\,i+2,\,j,\,j+1\}$; \\
        \>\>\> $i:=i+2$;
      \end{tabbing}
    \end{minipage}
  \end{compactenum}
  \caption{A Hamilton path $\pit_m$ on the graph of $\Qt_m$
    that induces an AOF-orientation $(n:=m+5)$.}
  \label{fig:path}
\end{figure}

We now verify that $\pit_m$ induces an AOF orientation on the graph
of~$\Qt_m$.  The $h$-vector of a simple polar-to-neighborly
$d$-dimensional polytope with $n=m+d+1$ facets is given by
$h_k={n-d-1+k\choose k} = {m+k\choose k}$ for
$k=0,1,\dots,d$. Therefore, by Proposition~\ref{prop:q_m:vertices},
\[
    \big(|T^1_m|,\,|T^3_m|,\,|T^4_m|,\,|T^2_m|,\,|T^0_m|\big) \;=\;
    (h_0,\,h_1,\,h_2,\,h_3,\,h_4).
\]
By Proposition~\ref{prop:hk}, it suffices to verify using
Figure~\ref{fig:q49-graph+path} that if the orientation of each edge
of the graph of~$\Qt_m$ is consistent with the total ordering induced
by~$\pit_m$, then the vertices of $T^1$, $T^3$, resp.~$T^4$ all have
in-degree $0$, $1$ resp.~$2$, furthermore $T^0$ and all but one of the
vertices of~$T^2$ have in-degree~$3$, and this vertex, the sink, has
in-degree~$4$.
\end{proof}

\section{Realizing the monotone Hamilton paths}
\label{sec:realizing-paths}

In this section we prove Theorem~\ref{thm:detailed}, and therefore
our Main Theorem.

\subsection{Outline of the inductive construction}
\label{subsec:outline}

For all $m\ge 0$, we first find an oriented hyperplane~$H_{m+1}$ that
separates the \emph{odd part} $T^{\le4}_{\rm odd}(m)=T_m^1\cup
T_m^3$ from the \emph{even part} $T^{\le4}_{\rm even}(m)=T_m^0\cup
T_m^2\cup T_m^4$ of $\pi_m$. We then create an
intermediate pair $(Q_{m+1}',\F_{m+1}')$ as in
Proposition~\ref{prop:q_m:vertices}: $Q_{m+1}':=Q_m\cap
H_{m+1}\pos$~is a simple polar-to-neighborly polytope of the same
combinatorial type as~$\Qt_{m+1}$, and the flag~$\F_{m+1}'$ of faces
is defined as in Definition \ref{def:Qt}(b). 

Our combinatorial model~$\Qt_{m+1}$ provides us with a Hamilton
path~$\pi_{m+1}$ on~$Q_{m+1}'$ that is not yet monotone with respect
to the objective function $f:\x\mapsto x_4$.  However, we will
choose~$H_{m+1}$ in such a way that there exists a pencil
\[
    \cH \ = \ \big\{H_t:t\in\Pr^1(\R)\cong\R\cup\{\infty\}\big\}
\]
of hyperplanes in~$\R^4$ with the following properties:

\begin{enumerate}[\quad({S}1)]

\item The common intersection of all hyperplanes in~$\cH$ is a
$2$-flat 
$
   R = \bigcap_{t\in \Pr^1(\R)}H_t
$
(the \emph{axis} of~$\cH$), and $\vertices Q_{m+1}'\cap R=\emptyset$.

\item The pencil $\cH$ ``sorts the vertices of~$Q_{m+1}'$ correctly'':
If $p\in H_r$ and $q\in H_s$ are vertices of~$Q_{m+1}'$ with
$r,s\ne\infty$ and $p$~precedes~$q$ in~$\pi_{m+1}$, then~$r<s$.
\end{enumerate}

We then apply a projective transformation~$\psi$ to
$\R^4\subset\Pr^4(\R)$ that maps $H_\infty$ to the hyperplane at
infinity. Because the common intersection~$R$ of all hyperplanes
in~$\cH$ is also mapped to infinity,
the image $\psi(\cH^b)=\psi(\cH\setminus
H_\infty)=\{\psi(H_t):t\in\R\}$ is a family of parallel affine
hyperplanes in~$\R^4$. The new objective function~$f$ is then defined
by the common normal vector to the hyperplanes in~$\psi(\cH^b)$, and the
Hamilton path~$\psi(\pi_{m+1})$ on $Q_{m+1}:=\psi(Q'_{m+1})$ is
strictly monotone with respect to~$f_{m+1}$ by~(S2).

\subsection{Properties of the family of polytopes}

\begin{notation} \label{conv:vertices}
  We use the following names for some special vertices of~$\Qt_m$:
  \begin{compactitem}[\quad$\triangleright$]
  \item The source $\{n-3,\,n-2,\,n-1,\,n\}$ of $\pit_m$ 
        is called~$\alpha_m$ (so that $T^1_m=\{\alpha_m\}$). 
  \item The sink is $\omega_m:=\{n-5,n-3,n-1,n\}\in
  T^2_m$. 
  \item  $\beta_m:=\{n-4,n-2,n-1,n\}$ (so that $T^0_m=\{\beta_m\}$).
  \item $\tau_m:=\{n-5,n-3,n-2,n-1\}\in T^4_m$.
  \end{compactitem}
\end{notation}

\begin{proposition}\label{prop:combinatorics} \ \samepage
  \begin{compactenum}[(a)]
  \item\label{prop:combinatorics:a} The induced subgraph of
    $sk^1(Q_m)$ on $T^1_m\cup T^3_m$ is a path of length~$m+1$ on the
    $m+2$ vertices $v_0^m=\alpha_m$, $v_1^m$, \dots, $v_{m+1}^m$, and the
    induced subgraph on~$T^2_m$ is a path
    $w_1^m,w_2^m,\dots,w_{m+1}^m$.

  \item\label{prop:combinatorics:b} For $0\le i \le m$, the
    edge~$e_i=\conv\{v_i^m,v_{i+1}^m\}$ in~$T^3_m$ is incident to a
    $2$-face~$G_i$ of~$Q_m$ such that the vertices of~$G_i\setminus
    e_i$ are consecutive in~$\pi_m\cap T^4_m$.

  \item\label{prop:combinatorics:c} For $1\le i \le m$, the edge~$f_i$
    of $Q_m$ that connects $w^m_i$~and~$w^m_{i+1}$ in~$T^2_m\cap\pi_m$ is
    incident to a quadrilateral~$R_i$ whose other two vertices are
    consecutive in~$T^4_m\cap\pi_m$.

  \item\label{prop:combinatorics:d} Set $G(m)=\vertices\bigcup_{i=0}^m
    G_i\setminus e_i$ and $R(m)=\vertices\bigcup_{i=1}^m R_i\setminus f_i$. Then
    $G(m)\cup R(m)=T^4_m$, and $G(m)\cap R(m)=\tau_m$.
    
  \end{compactenum}
\end{proposition}

\begin{figure}[htbp]
  \centering
  \input{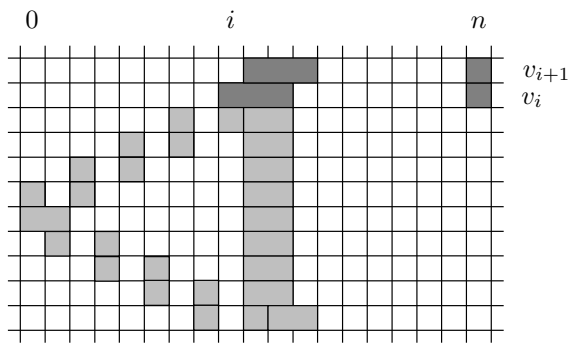}
  \caption{Vertices of a $2$-face incident to $v_i=\{i,i+1,i+2,n\}$
    and $v_{i+1}=\{i+1,i+2,i+3,n\}$ (dark) in~$T^1_m\cup T^3_m$. The
    light vertices lie in $T^4_m$ and form a subpath of~$\pi_m$.}
  \label{fig:lemma-combinatorics}
\end{figure}

\begin{figure}[htbp]
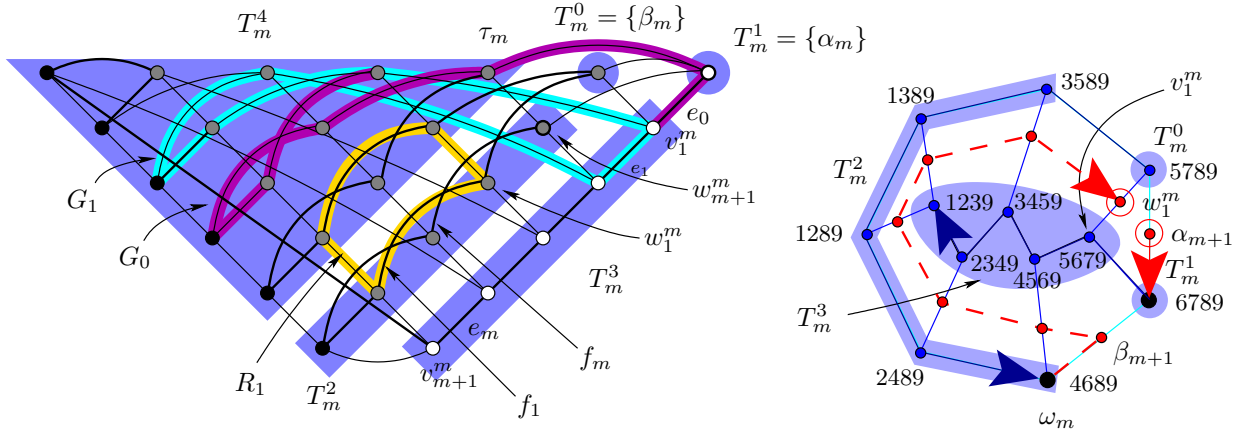

  \centering
  \begin{minipage}{.64\linewidth}
     \input{q49-graph+path+2-faces.pstex_t}
  \end{minipage}\hfill
  \begin{minipage}{.35\linewidth}
     \raisebox{-6cm}{\input{facet-8-cut.pstex_t}}
  \end{minipage}
  \caption{\emph{Left:} More details about the graph of~$Q_m$. We have
    highlighted the graphs of the $2$-faces $G_0$~and~$G_1$ that
    correspond to the edges $e_0$~and~$e_1$ by
    Proposition~\ref{prop:combinatorics}~\eqref{prop:combinatorics:b},
    and the $2$-face~$R_1$ that corresponds to the edge~$f_1$
    according to
    Proposition~\ref{prop:combinatorics}~\eqref{prop:combinatorics:c}.
    \emph{Right:} The portion of the new Hamilton path~$\pit_{m+1}$
    in the facet~$F^3_m$. } \label{fig:graph+path+2-faces}
\end{figure}

\begin{proof} \eqref{prop:combinatorics:a} All vertices of
$T^3_m$ are of the form $\{i,i+1,i+2,n\}$ for~$1\le i\le n-3$, and the
only way for two such vertices $v^m_i$~and~$v^m_j$ to be adjacent
for~$i<j$ is to have~$j=i+1$. The statement about the $w^m_i$ follows
in a similar way.  \eqref{prop:combinatorics:b} For $1\le i \le m+1$,
the $2$-face incident to $v_{m+2-i}=\{i,i+1,i+2,n\}$ and
$v_{m+1-i}=\{i+1,i+2,i+3,n\}$ that is the intersection of the facets
$i+1$ and~$i+2$ consists of the vertices of
Figure~\ref{fig:lemma-combinatorics}.  The claim
\eqref{prop:combinatorics:b} follows because these vertices form a
contiguous segment of~$\pi_m$, and
\eqref{prop:combinatorics:c}~and~\eqref{prop:combinatorics:d} from
Figure~\ref{fig:graph+path+2-faces} (left). 
\end{proof}

\begin{observation}
  The new start vertex $\alpha_{m+1}$ of $\pit_{m+1}$ lies on
  $\conv\{\alpha_m,\beta_m\}$, the new end vertex~$\omega_{m+1}$
  on~$\conv\{v_1^m,\beta_m\}$, and $\beta_{m+1}$
  on~$\conv\{\alpha_m,\omega_m\}$; see
  Figure~\ref{fig:graph+path+2-faces} (right).
\end{observation}

\subsection{Start of the induction and inductive invariant}

We work in $\R^4$ with standard coordinate vectors $e_1,e_2,e_3,e_4$.
An essential tool will be \emph{shear transformations}: these are
linear maps $\sigma_{i,j}^a:\R^4\to\R^4$ for $i,j\in\{1,2,3,4\}$,
$i\ne j$, and~$a\in\R$ whose matrix is $I_4+a\delta_{i,j}$ with
respect to the standard basis of~$\R^4$. Here $I_4$~is the $4\times4$
unit matrix and $\delta_{i,j}$~is the $4\times 4$ matrix whose only
nonzero entry is a $1$ in position~$(i,j)$. In particular,
$\sigma_{i,j}^a$ maps $e_i$ to $e_i+ae_j$, and the standard basis
vectors $e_k$, $k\ne i$, to themselves.

The start of the induction is the pair $(Q_0,\F_0)$, where $Q_0$ is
the $4$-simplex whose vertices $v_1,v_2,v_3,v_4,v_5$ are given by the
columns of the matrix
\begin{equation}\label{eq:q0}
  \begin{pmatrix}
    0 & 0   &   1 & 0 & 0\\
    0 & 1   &   0 & 0 & 0\\
    -3 & -1 &   3 & 2 & 1\\
    -2 & -1/2 & 0 & 1/4 & 2
  \end{pmatrix}\,,
\end{equation}
and $\F$~is the  flag $\F_0:F_0^0\subset
F_0^1\subset\dots\subset F_0^4=Q^4_0$ of faces labeled as in
Definition~\ref{def:Qt}. In particular, the vertices~$v_i$ lie in
the following tips,
\begin{center}\small
  \begin{tabular}[c]{c|c|c|c|c}
    $v_1$ & $v_2$ & $v_3$ & $v_4$ & $v_5$\\\hline
    $T_0^1$ & $T^3_0$ & $T^4_0$ & $T^0_0$ & $T^2_0$
  \end{tabular}\;,
\end{center}
$F_0^2=\conv\{v_1,v_4,v_5\}$, $F_0^3=\conv\{v_1,v_2,v_4,v_5\}$, and
 $\pi_0=(v_1,v_2,v_3,v_4,v_5)$.  

\smallskip
For all $m\ge0$ the polytopes $Q_m$ will maintain the following
property:

\begin{compactenum}[\quad({M}1)]
\item\label{pre-ascending} The Hamilton path $\pi_m$ in the
  $1$-skeleton of~$Q_m$ is strictly monotone with respect to the
  objective function~$f:\R^4\to\R$, $\x\mapsto x_4$.
\end{compactenum}

\subsection{Induction step I: Positioning the polytope}
\label{subsec:step1}

In this and the following section, we will position the polytope~$Q_m$
in such a way that the coordinate subspaces of~$\R^4$ are compatible
with the flag~$\F_m$. More precisely, 
\begin{itemize}[\quad$\triangleright$]
\item $F^3_m = Q_m\cap\{\x\in\R^4:x_1=0\}$, and
  $T^4_m\subset\{\x\in\R^4:x_1>0\}$; and
\item the hyperplane $H_S=\{\x\in\R^4:x_3=0\}$ will separate $T^{\le
4}_{\rm even}(m)$ from $T^{\le 4}_{\rm odd}(m)$.
\end{itemize}

\begin{lemma}\label{lem:position-F2}
  Let $\pi$ be the linear projection $\pi:\R^4\to\ _\R\langle e_3,e_4\rangle$,
  and use the notation of Convention~\ref{conv:vertices} and
  Proposition~\ref{prop:combinatorics}\eqref{prop:combinatorics:a}.
  Then there exists a non-singular affine transformation~$\sigma$ of~$\R^4$
  such that $Q_m\equiv\sigma(Q_m)$~satisfies the following
  additional conditions, while $\pi_m\equiv\sigma(\pi_m)$~still
  satisfies~\plabel{\ref{pre-ascending}}: 

  \begin{enumerate}[\quad({M}1)]
    \setcounter{enumi}{1}
  \item\label{pre-f2} $F^2_m\subset\{\x\in\R^4:x_1=0\}$.

  \item\label{pre-f3} $\aff F^3_m=\{\x\in\R^4:x_1=0\}$ and
    $Q_m\subset\{\x\in\R^4: x_1\ge0\}$.

  \item\label{pre-f4} $(\alpha_m)_2=0$, $r_2<0$ for all $r\in
    F^2_m\setminus\{\alpha_m\}$, and $(\beta_m)_2<(v^m_1)_2$. 

  \item \label{ind:post-1} The image of~$F^2_m$ under~$\pi$ is
    full-dimensional: $\dim\aff\big(\pi(F^2_m)\big)=2$.

  \item \label{ind:post-4} The $3$-flat $H_S=\{\x\in\R^4:x_3=0\}$
    strictly separates $T^{\le 4}_{\rm even}(m)$ from $T^{\le 4}_{\rm
    odd}(m)$.  Moreover, we may choose the point of~$H_S\cap F^3_m$ of
    lowest $4$-coordinate to be
    $\alpha_{m+1}=\conv\{\alpha_m,\beta_m\}\cap H_S$,
    where~$(\alpha_{m+1})_4=\tau_4:=(\tau_m)_4$.  
  \end{enumerate}
\end{lemma}

\noindent\emph{Proof.} Properties \plabel{\ref{pre-f2}}
and~\plabel{\ref{pre-f3}} are a matter of trivial affine transforms
that can be chosen to leave the $4$-coordinates invariant,
thereby maintaining~\plabel{\ref{pre-ascending}}, and
property~\plabel{\ref{pre-f4}} can be achieved via a translation and a
shear~$\sigma^a_{2,4}:x_2\mapsto x_2+ax_4$.

For~\plabel{\ref{ind:post-1}}, choose $t\in F^2_m$ with $t_4=q_4$ for
some $q\in T^3_m$; such a point exists, since $\alpha_m\in F^2_m$, and
$(\alpha_m)_4 < q'_4 < \max\{s_4:s\in F^2_m\}$ for all $q'\in T^3_m$
by~\plabel{\ref{pre-ascending}} and Remark~\ref{rem:order-pi}.
Translate~$t$ 
such that $t=(0,t_2,0,0)$ with
$t_2<0$, and apply a shear transform $\sigma_{3,2}^b:x_3\mapsto
x_3+bx_2$ to~$\R^4$, where $b\in\R$~is chosen such that
$\pi\big(\sigma_{3,2}^b(q)\big)= \pi\big(\sigma_{3,2}^b(t)\big)$.
This can be done because
$\pi(t)-\pi(q)\in\R\pi(e_3)$.  Then \plabel{\ref{ind:post-1}}~is
fulfilled because $\dim\aff F^3_m=3$: supposing that
$\dim\aff\big(\pi(F^2_m)\big)=1$ would imply via~$t\in F^2_m$
and~$q\in F^3_m$ that~$q\in\aff F^2_m$; however, this is absurd by the
choice~$q\in T^3_m$.  Note that none of the maps we used
affects~\plabel{\ref{pre-f2}}--\plabel{\ref{pre-f4}}.

\begin{figure}[htbp]
  \centering
  \psfrag{am1}{$\alpha_{m+1}$} \psfrag{alpha}{$\alpha_m=T^1_m$} 
  \psfrag{wm}{$\omega_m$} \psfrag{betam}{$\beta_m$}
  \psfrag{om1}{$\omega'$}
  \psfrag{ell}{$\ell$} \psfrag{ell1}{$\ell'=\sigma(\ell)$}
  \psfrag{T0}{$\beta_m=T^0_m$} \psfrag{T3}{$T^3_m$}
  \psfrag{F2}{$F^2_m$}
  \psfrag{sigma}{$\sigma^b_{3,4}$}
  \psfrag{bm1}{$\beta_{m+1}$} \psfrag{3}{$e_3$} \psfrag{4}{$e_4$}
  \includegraphics[width=.9\linewidth]{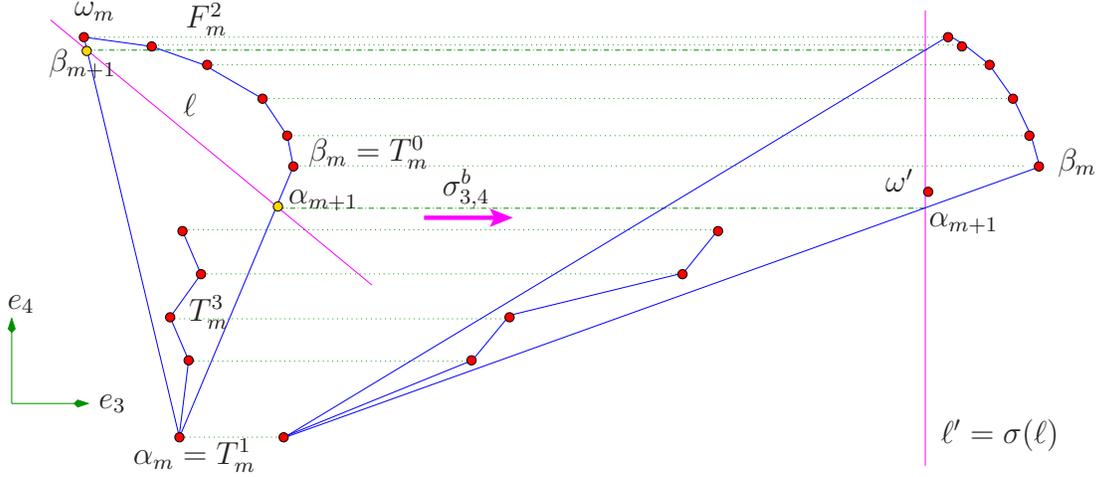} 
  \caption{Positioning the polytope, step~\plabel{\ref{ind:post-4}}. 
    The map~$\sigma^b_{3,4}$ shears the polytope until (the preimage
    under~$\pi$ of) a vertical line~$\ell'$ separates the odd from the
    even tips. On the right, the approximate position of~$\omega'$ is
    marked; cf.~Lemma~\ref{lem:asymptotics}.}
    \label{fig:positioning}
\end{figure}

For~\plabel{\ref{ind:post-4}}, define $\bt$ to be the point of
greatest $3$-coordinate of~$F^2_m\cap\{\x\in\R^4:x_4=\tau_4\}$. In
particular, $\bt_4>\max_{z\in T^3}z_4$ by~\plabel{\ref{pre-ascending}},
and $\bt$~lies either on the edge~$\conv\{\alpha_m,\beta_m\}$ or on
the edge $\conv\{\alpha_m, \omega_m\}$ of~$F^2_m\subset Q_m$
(cf.~Figure~\ref{fig:positioning}).

\noindent Possibly using the transform $x_3\mapsto -x_3$, we can
achieve $\bt\in\conv\{\alpha_m,\beta_m\}$, and $\bt=\alpha_{m+1}$
after a translation along the $3$-axis.  Now choose a non-horizontal
line~$\ell$ through~$\alpha_{m+1}$ such that $\pi(\ell)$ separates
$\pi(T^1_m\cup T^3_m)$ from $\pi(F^2_m\setminus T^1_m)$ (for example,
perturb $\ell=\alpha_{m+1}+\R e_3$), translate~$Q_m$ again such that
$\alpha_{m+1}=0$, and apply a shear $\sigma_{3,4}^c:x_3\mapsto
x_3+cx_4$ to~$\R^4$ such that
$\ell':=\sigma_{3,4}^c(\ell)=\{\x\in\R^4:x_1=x_3=0\}\cap\aff F^2_m$ is
vertical, and $x_3<0<y_3$ for all $x\in T^1_m\cup T^3_m$ and $y\in
F^2_m\setminus T^1_m$ (cf.~Figure~\ref{fig:positioning}).  If the
hyperplane $\pi^{-1}\big(\pi(\ell')\big)$ does not yet
separate~$T^1_m$ from~$T^4_m$, apply another shear
$\sigma_{3,1}^d:x_3\mapsto x_3+ d x_1$ with~$d>0$ until it does (note
that~\plabel{\ref{pre-f3}} already holds), and then define
$H_S:=\pi^{-1}\big(\pi(\ell')\big)$.  This hyperplane then separates
the odd and even parts of~$\pi_m$ by construction,
and~$(\alpha_{m+1})_4=\tau_4$ also by construction and because the
shears $\sigma_{3,4}^c$~and~$\sigma_{3,1}^d$ do not affect
$4$-coordinates.  Neither do they affect
conditions~\plabel{\ref{pre-ascending}}--\plabel{\ref{ind:post-1}}, so
we define~$\sigma$ as the composition of all these
maps. \hfill $\Box$

\begin{remark}
  The
  conditions~\plabel{\ref{pre-ascending}}--\plabel{\ref{ind:post-4}}
  are satisfied by the coordinates~\eqref{eq:q0} for~$Q_0$.
\end{remark}

\subsection{Induction step II: Finding the cutting plane}
\label{subsec:step2}

In this section, we will find a hyperplane~$H_{m+1}$ that gives rise
to a polytope~$Q'_{m+1}=Q_m\cap H_{m+1}\pos$ of the same combinatorial
type as~$\Qt_{m+1}$. Namely, assume that
\plabel{\ref{pre-ascending}}--\plabel{\ref{ind:post-4}} hold, define
$H_{m+1}$ to be the hyperplane $\{\x\in\R^4:\n^T\!\x=0\}$ with
$\n=(0,-\delta,1,\e)^T$ for some small $\e\gg\delta>0$, and assign the
label~$n+1=m+d+2$ to~$H_{m+1}$.  Note that $H_{m+1}$ converges
to~$H_S$ as~$\e,\delta\to 0$.

\begin{remark}
  Up to now, we have put the facet~$F_m^3$ into the $3$-plane
  $\{\x\in\R^4:x_1=0\}$ and the tip~$T^4$ into the half-space
  $\{\x\in\R^4:x_1>0\}$. This allows us to move ``almost all'' of
  the vertices of $\pi_m$ (namely, the portion inside~$T^4_m$) ``out
  of the way'', via a shear~$\sigma_{3,1}^a$ that only affects
  $3$-coordinates. These ``old'' vertices will be dealt with in
  Lemma~\ref{lem:sweep} below.
  
  We still need to arrange for the first and last part of~$\pi_{m+1}$
  to be traversed in the right order. We achieve this by adjusting the
  position of~$H_{m+1}$ via the parameters $\e$~and~$\delta$ in the
  definition of~$\n$ (note that we chose $n_1=0$, because we are
  already done with~$T^4_m$). If~$\delta=0$, then $\pi(H_{m+1})$ is a
  line whose slope is determined by~$\e$. We choose $\e>0$ to `push
  out' the first part $T^1_{m+1}\cup T^3_{m+1}$ of the new
  path~$\pi_{m+1}$.  However, if we left $\delta=0$ we would not
  correctly sweep the last portion $T^0_{m+1}\cup T^2_{m+1}$.
  Items~\plabel{\ref{lem:asymptotics:4}}--\plabel{\ref{lem:asymptotics:6}}
  of Lemma~\ref{lem:asymptotics} guarantee a correct sweep in
  Lemma~\ref{lem:sweep} for sufficiently small $0<\delta\ll\e$.
\end{remark}

\begin{lemma}\label{lem:asymptotics}
  Assume
  conditions~\plabel{\ref{pre-ascending}}--\plabel{\ref{ind:post-4}}
  and $(\alpha_{m+1})_3=(\alpha_{m+1})_4=0$, and fix vertices $q\in
  T^{\le4}_{\rm odd}(m)$ and $s\in T^{\le4}_{\rm even}(m)$. Let
  $q'=\conv\{q,s\}\cap H_{m+1}$ be the intersection with~$H_{m+1}$ of
  the line through $q$~and~$s$ (which is not necessarily an edge
  of~$Q_m$).  Then, if $a>0$~is sufficiently large and
  $0<\delta\ll\e$~are sufficiently small, the image
  $\sigma_{3,1}^a(Q_m)$ of~$Q_m$ under the shear $\sigma_{3,1}^a$
  satisfies the following
  conditions~\plabel{\ref{lem:asymptotics:2}}--\plabel{\ref{lem:asymptotics:6}};
   cf.~also Figure~\ref{fig:inductive-step} below.

  \begin{enumerate}[\quad({M}1)] \setcounter{enumi}{6}

  \item\label{lem:asymptotics:2} $q_3'>0$ for
    $0<\delta\ll\e$, and $q_3'\searrow 0$ as $\delta,\e\searrow0$. In other
    words, all points in $\sigma_{3,1}^a(Q_m)\cap H_{m+1}$ can be
    chosen to have positive $3$-coordinate, but to lie arbitrarily
    close to~$H_S$.

  \item\label{lem:asymptotics:4} Set
    $u:=(v_1^{m+1})'=\conv\{\alpha_m,\tau_m\}\cap H_{m+1}$ and suppose
    that $q'\ne u$. Then the
    image $\pi\big(\!\aff\{u,q'\}\big)\subset\ _\R\langle e_3,e_4\rangle$ of the
    line through $u$~and~$q'$ under~$\pi$ comes arbitrarily close to
    being vertical as $a\to\infty$ and $\e,\delta\to0$.

    \item\label{lem:asymptotics:5} Set
      $\alpha':=\alpha'_{m+1}=\conv\{\alpha_m, \beta_m\}\cap
      H_{m+1}$. If $q,\bar q\in T^3_m$ and $q_4<\bar q_4$, so that
      $q',\bar q'\in T^3_{m+1}$ and $q'_4 < \bar q'_4$, then the slope
      $\sigma_{\alpha' \bar q'}$ of the
      line~$\pi\big(\!\aff\{\alpha', \bar q'\}\big)$ is greater than
      the slope $\sigma_{\alpha' q'}$ of the
      line~$\pi\big(\!\aff\{\alpha',q'\}\big)$ (and both are
      negative).

      \item\label{lem:asymptotics:6} Set
      $\omega':=\omega'_{m+1}=\conv\{\beta_m, v_1^m\}\cap
      H_{m+1}$. Then the slope $\sigma_{\omega'\alpha'}$
      of~$\pi\big(\!\aff\{\omega',\alpha'\}\big)$ is less than the
      slope $\sigma_{\omega' u}$ of~$\pi\big(\!\aff\{\omega',
      u\}\big)$.  

  \end{enumerate}
\end{lemma}

\begin{proof} We abbreviate
  $\sigma=\sigma_{3,1}^a$.  For~\plabel{\ref{lem:asymptotics:2}}, we
  have $\conv\{q,s\}\cap H_{m+1}\ne\emptyset$ since $q$~and~$s$ are
  separated by~$H_{m+1}$ for small enough $\delta,\e$.  We calculate
  the intersection point~$q'=\conv\{q,s\}\cap H_{m+1}$ by solving
  $\n^T\!\q + \mu\,\n^T\!(\s-\q) = 0$ for~$\mu$, obtaining
  \[
     \q' \ = \  \q + \frac{\n^T\!\q}{\n^T\!(\q-\s)}\,(\s-\q).
  \]
  By~\plabel{\ref{pre-f2}}, the map $\sigma$ leaves the points
  $\alpha'$, $q$, and~$\omega'$ invariant, and maps $\s$~to
  $\sigma(\s)=\s+as_1\be_3$; as a consequence,
  $\n^T\!\sigma(\s)=\n^T\!\s+as_1$.  Using $\n^T\!\q = -\delta q_2 +
  q_3+\e q_4$, we obtain
\begin{eqnarray}
  \sigma(\q') &=&
  \q+\frac{\n^T\!\q}{\n^T(\q-\s)-as_1}\,(\s-\q+as_1\be_3)\notag\\[1ex]
   &\xrightarrow[\;a\to\infty\;]{}&
   \q+ (0,\,0,\,-\n^T\!\q,\,0)^T
   \ = \  (0,\, q_2,\, \delta q_2 -\e q_4,\, q_4)^T. \label{eq:star1}
\end{eqnarray}
Because $q_4<(\alpha_{m+1})_4=0$, we can choose $0<\delta\ll\e$ so
small that $\sigma(\q')_3>0$ (note that~$q_2\le 0$
by~\plabel{\ref{pre-f4}}). In particular, we obtain
$\sigma(\q')_3\searrow 0$ as~$\e,\delta\searrow0$.

\smallskip
Statement \plabel{\ref{lem:asymptotics:4}} follows from~\eqref{eq:star1} and
the fact that
\[
    \lim_{a\to\infty}
    \,\frac{\sigma(\q')_4-\sigma(\bu)_4}{\sigma(\q')_3-\sigma(\bu)_3}
    \ = \ \frac{q_4-u_4}{\delta(q_2-u_2) - \e(q_4-u_4)}\;.
\]

\smallskip
For~\plabel{\ref{lem:asymptotics:5}}, note that since $\alpha'$ is
invariant under $\sigma$,
\[
    \sigma_{\alpha' q'}  \ = \
    \frac{\sigma(\q')_4-\alpha'_4}{\sigma(\q')_3-\alpha'_3} 
    \quad\xrightarrow[\;a\to\infty\;]{}\quad
    \frac{q_4-\alpha'_4}{\delta q_2 - \alpha'_3 - \e q_4}\,, 
\]
and similarly for $\bar q$; the statement now follows from $q_4<\bar
q_4$ and $0<\delta\ll\e$. 

\smallskip
To prove~\plabel{\ref{lem:asymptotics:6}}, set $\alpha:=\alpha_m$,
$\beta:=\beta_m$, $v:=v_1^m$ and $\tau:=\tau_m$. Then
$u=\conv\{\alpha,\tau\}\cap H_{m+1}$,
$\alpha'=\conv\{\alpha,\beta\}\cap H_{m+1}$, and
$\omega'=\conv\{v,\beta\}\cap H_{m+1}$. We need to verify that
\[
    \sigma_{\omega'\alpha'} \ :=\  
    \frac{\alpha'_4-\omega'_4}{\alpha'_3-\omega'_3} \ <\ 
    \frac{u_4-\omega'_4}{u_3-\omega'_3} \ =:\  
    \sigma_{\omega' u}.
\]
From equation~\eqref{eq:star1} and condition~\plabel{\ref{pre-f4}}, we deduce
that $\lim_{a\to\infty} u = (0, 0, -\e \alpha_4, \alpha_4)^T$.  For
$\alpha'$ and~$\omega'$ we get the following expressions:
  \begin{gather*}
    \balpha' \ =\ 
    \balpha+\frac{\n^T\!\balpha}{\n^T\!(\balpha-\bbeta)}\,(\bbeta-\balpha) \ =\ 
    \begin{pmatrix}
      0\\ 0\\ \alpha_3 \\ \alpha_4
    \end{pmatrix}
    + \frac{\alpha_3 + \e\alpha_4}
           {\delta\beta_2 + \alpha_3-\beta_3 + \e(\alpha_4-\beta_4)}
    \begin{pmatrix}
      0\\ \beta_2 \\ \beta_3-\alpha_3  \\ \beta_4-\alpha_4 
    \end{pmatrix},
    \\
    \bomega' \ =\ 
    \bv+\frac{\n^T\!\bv}{\n^T\!(\bv-\bbeta)}\,(\bbeta-\bv) \ =\ 
    \begin{pmatrix}      
      0\\ v_2 \\ v_3 \\ v_4
    \end{pmatrix}
    + \frac{-\delta v_2 + v_3 +\e v_4}
           {-\delta(v_2-\beta_2) + v_3-\beta_3 + \e(v_4-\beta_4)}
    \begin{pmatrix}
      0\\ \beta_2-v_2 \\ \beta_3-v_3  \\ \beta_4-v_4 
    \end{pmatrix}.
  \end{gather*}
For convenience, we will verify that
$1/\sigma_{\omega'\alpha'}>1/\sigma_{\omega' u}$. Indeed, expanding
these expressions in terms of $\delta,\e$, we obtain 
\begin{eqnarray*}
  \frac{1}{\sigma_{\omega'\alpha'}} &=&
    \frac{\beta_3v_2 - \beta_2v_3 +
      \smash{\overbrace{\alpha_3(\beta_2-v_2)}^{t_1}}\mathstrut}
         {v_3(\alpha_4-\beta_4) + \beta_3(v_4-\alpha_4) + 
           \underbrace{\alpha_3(\beta_4-v_4)}_{t_2}} \,\delta 
         -\e + p_1(\delta,\e),\\
  \frac{1}{\sigma_{\omega' u}} &=&
    \frac{\beta_3v_2 - \beta_2v_3}{v_3(\alpha_4-\beta_4) +
      \beta_3(v_4-\alpha_4)} \,\delta - \e + p_2(\delta,\e),
\end{eqnarray*}
where $p_1$ and $p_2$ are power series in $\delta,\e$ with min-degree
at least~$2$.  Notice that up to terms of degree at least~$2$ in
$\delta,\e$, the two formulas are equal except for the
expressions~$t_1$ resp.~$t_2$ in the numerator resp.\ denominator
of~$1/\sigma_{\omega'\alpha'}$. Therefore, we can write the difference
between the inverses of the slopes as
\[
  \frac{1}{\sigma_{\omega'\alpha'}} -  \frac{1}{\sigma_{\omega' u}} \ =\ 
  \left(\frac{A+t_1}{B+t_2} - \frac{A}{B}\right)\delta + p_3(\delta,\e).
\]
Since $\alpha_3<(\alpha_{m+1})_3<0$ by assumption and $\beta_2<v_2$
by~\plabel{\ref{pre-f4}}, we obtain $t_1>0$; and the inductive
assumption~\plabel{\ref{pre-ascending}} implies that $\beta_4>v_4$ and
therefore $t_2<0$. The claim follows.
\end{proof}

\subsection{Induction step III: The projective transformation}
\label{subsec:step3}

Finally, we construct a $1$-parameter family
$\cH=\{H_t:t\in\Pr^1(\R)\}$ of hyperplanes that contains a
$2$-plane~$R$ as their common ``axis'', as in
Section~\ref{subsec:outline}.
Let $O=\pi\big(b+\e_1(\omega-\alpha)-\e_3e_3\big)$ for some
small~$\e_1,\e_3>0$, so that $O$~lies outside but very close to the
edge $\conv\{\alpha,\omega\}$ of~$\pi(F^2_{m+1})$, and define the
$2$-plane $R\subset\R^4$ to be~$R=\pi^{-1}(O)$.

\begin{lemma}  \label{lem:sweep}
  Let $\cH$ be the pencil of hyperplanes in~$\R^4$ sharing the
  2-plane~$R$, and such that $\pi(H_\infty)$ is the line through~$O$
  parallel to~$\conv\{\alpha,\omega\}$, and the slope of $\pi(H_r)$ is
  smaller than the slope of~$\pi(H_s)$ exactly if~$r<s$.  Then
  $\cH$~fulfills~(S2), i.e., it sorts the vertices of~$Q_{m+1}$ in the
  order given by~$\pi_{m+1}$.
\end{lemma}

\begin{landscape}

\begin{figure}[htbp]
  \centering
  \psfrag{alpham}{$T^1_m=\{\alpha_m\}$} \psfrag{omegam}{$\omega_m$} 
  \psfrag{betam}{$T^0_m=\{\beta_m\}$}
  \psfrag{alpha}{$\alpha'$} \psfrag{beta}{$\beta$}
  \psfrag{omega}{$\omega'$} \psfrag{tau}{$\tau_m$}
  \psfrag{alpha1}{$\alpha'=T^1_{m+1}$} \psfrag{beta}{$\beta$}
  \psfrag{H}{$H_\infty$}
  \psfrag{u}{$u$} \psfrag{vm}{$v_m$}
  \psfrag{O}{$\pi(R)=O$}
  \psfrag{T2}{$T^2_{m+1}$} \psfrag{T3}{$T^3_{m+1}$}
  \psfrag{T3m}{$T^3_m$} \psfrag{F2m}{$F^2_m$} \psfrag{T4m}{$T^4_m$} 
  \psfrag{t4}{$\tau_4$}
  \psfrag{e3}{$e_3$} \psfrag{e4}{$e_4$} \psfrag{e12}{$e_1$, $e_2$}
  \psfrag{10}{\ref{lem:asymptotics:5}}
  \psfrag{11}{\ref{lem:asymptotics:6}}
  \includegraphics[width=19.2cm]{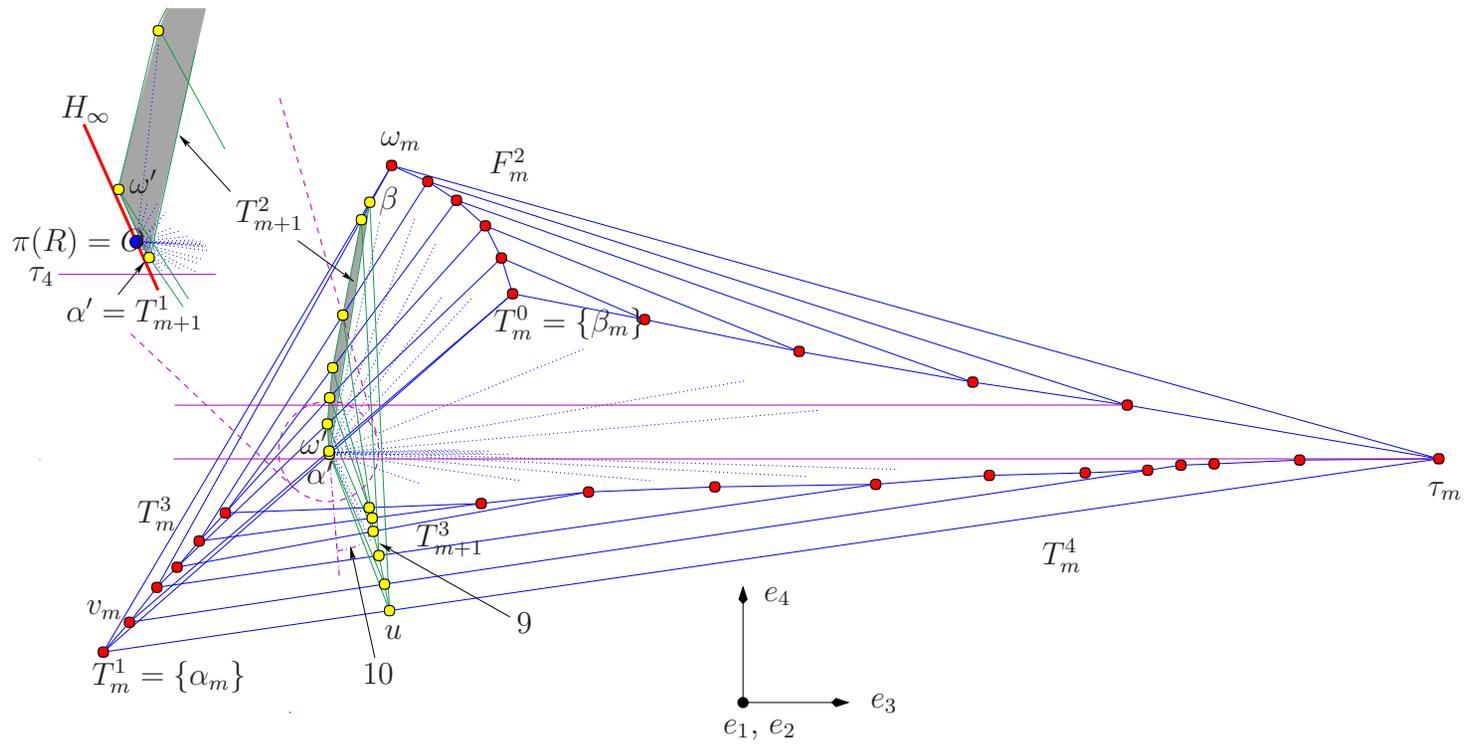}
  \caption{\emph{The inductive step:} We show the projection of the
    polytope~$Q_4$ to the~$\langle 3,4\rangle$-plane, and the vertices
    obtained by intersecting $Q_4$~with~$H_5$. The arrows next to the
    labels \ref{lem:asymptotics:5}~and~\ref{lem:asymptotics:6} point
    to the lines about whose slope the corresponding condition in
    Lemma~\ref{lem:asymptotics} makes an assertion.  The line through
    $O$ is the projection of the $3$-plane~$H_\infty$. A sweep
    around~$O$ encounters all vertices of $Q_m\cap H_{m+1}$ in the
    correct order~$\pi_m$ prescribed by~$\pit_{m+1}$.}
  \label{fig:inductive-step}
\end{figure}
\end{landscape}

\begin{proof} We examine the pieces of~$\pi_{m+1}$ in order;
cf.~Figure~\ref{fig:inductive-step}. 

\begin{compactitem}[$\triangleright$]
\item \emph{$T^1_{m+1}=\{\alpha\}$ is the start of~$\pi_{m+1}$:} This
  follows for small enough~$\e_3$ by~\plabel{\ref{lem:asymptotics:6}}.
  
\item \emph{$T^3_{m+1}$ is traversed next, in the right order, and
    before~$T^4_{m+1}$:} The first two statements follow
  from~\plabel{\ref{lem:asymptotics:2}},~\plabel{\ref{lem:asymptotics:4}}
  and~\plabel{\ref{lem:asymptotics:5}}, and the last one because $z_3\to\infty$
  as~$a\to\infty$ for any~$z\in T^4_{m}$, while the $3$-coordinates
  of~$T^3_{m+1}$ remain bounded by~\plabel{\ref{lem:asymptotics:2}}.

\item \emph{The correct order in~$T^4_m\subset T^4_{m+1}$.}  By
Proposition~\ref{prop:combinatorics}\eqref{prop:combinatorics:b}, each
of the edges $e_i=\conv\{v_i^m,v_{i+1}^m\}$, $0\le i \le m$, of~$T^1_m
\cup T^3_m$ is incident to an~$(m+1)$-gonal 2-face~$G_i$ (see
Figure~\ref{fig:graph+path+2-faces}), and the edges~$E_i$ of~$G_i$ not
incident to~$e_i$ form a monotone subpath of~$\pi_{m+1}$.  This
implies that for each $e_i\in T^3_m$, the slopes of the projection of
each~$E_i$ to~$\ _\R\langle e_3,e_4\rangle$ are strictly positive (and, by
convexity, monotonically decreasing; see Figure~\ref{fig:order-4a}).

  \begin{figure}[htbp]
    \centering\small
    \psfrag{alpha}{$\alpha$} \psfrag{tau}{$\tau$} \psfrag{O}{$O$}
    \psfrag{ei}{$e_i$} \psfrag{Gi}{$G_i$} \psfrag{E_i}{$E_i$}
    \includegraphics[width=.4\linewidth]{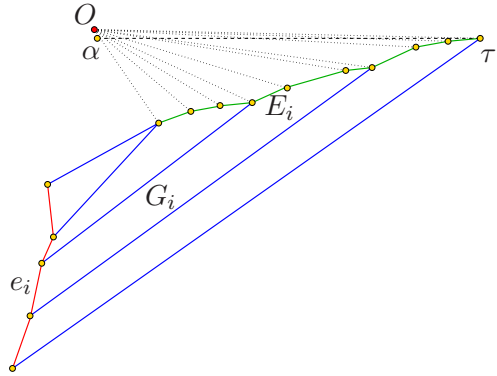}
    \caption{Convexity of the $(m+1)$-gonal faces enforces the correct
      order in~$T^4_{m}\subset T^4_{m+1}$.}
    \label{fig:order-4a}
  \end{figure}

Therefore, $\pi\big(\bigcup_{i=0}^m E_i\big)$ is a strictly increasing
chain of edges, and this remains true after applying the linear
map~$\sigma=\sigma_{3,1}^a$ by invariance of the $e_i$'s and all
$4$-coordinates under~$\sigma$, and the convexity of the projections
of $2$-faces.  The correct order up to~$\tau$ in~$T^4_m\subset
T^4_{m+1}$ follows from condition~\plabel{\ref{ind:post-4}}:
$\alpha_4\ge s_4$ for all $s\in\bigcup_{i=0}^m \vertices
G_i\setminus\vertices e_i$.  Similarly, the $4$-gonal $2$-faces
incident to~$T^2_m$ of
Proposition~\ref{prop:combinatorics}\eqref{prop:combinatorics:c}
enforce the right order between~$\tau$ and~$T^0_m$.
  
\item \emph{$T^2_{m+1}$ is traversed after~$T^4_{m+1}$:}
  Since~$\beta$, the first vertex of~$\pi_{m+1}$ to come
  after~$T^4_{m+1}$, lies on $\conv\{\alpha_m,\omega_m\}$, this can be
  achieved by choosing $\e$~and~$\e_1$ suitably small.
  
\item \emph{Correct order in $T^2_{m+1}$ and~$T^0_{m+1}$.}  This
  follows because the convex polygon $\pi(F^2_{m+1})$ is star-shaped
  with respect to any point on its boundary, and the choice of~$O$
  close to an edge of~$\pi(F^2_{m+1})$.

\end{compactitem}

\noindent This concludes the proof of Lemma~\ref{lem:sweep}.
\end{proof}

\smallskip

Finally, we apply the projective transform~$\psi:\R^4\to\R^4$,
$\x\mapsto\x/(\a\x-a_0)$ that sends the $3$-plane
$H_\infty=\{\x\in\R^4:\a\x=a_0\}$ to infinity, and set
$Q_{m+1}:=\psi(Q'_{m+1})$. Lemma~\ref{lem:sweep} then implies the
inductive condition~\plabel{\ref{pre-ascending}}, namely that
$Q_{m+1}$~admits an monotone Hamilton path~$\pi_{m+1}$.  The proof
of Theorem~\ref{thm:detailed}, and so of the Main Theorem, is
concluded. \hfill $\Box$

\section{Acknowledgements}

It is a pleasure to thank G\"unter M. Ziegler for suggesting this
problem, and Volker Kaibel for his careful reading of an earlier
version of the paper.


\end{document}